# GROUPOIDS IN DIFFEOLOGY

PATRICK IGLESIAS-ZEMMOUR

Abstract. This expository paper recounts the development and application of the concept of the diffeological groupoid, from its introduction in 1985 to its use in current research. We demonstrate how this single concept has served as a powerful and unifying tool for defining fundamental structures, analyzing the stratification of complex spaces like orbifolds, building a bridge to noncommutative geometry, and, most recently, forging new approaches to geometric quantization. The paper aims to provide a cohesive narrative of this journey, making explicit certain concepts like the "Klein groupoid" and showcasing the enduring vitality of the diffeological groupoid in modern geometry and physics.

## Introduction

At first glance, a paper focused on *diffeological* groupoids might seem to sit adjacent to the theme of topological group theory. However, this paper will argue that the opposite is true: diffeology provides a powerful smooth enrichment of topology, and its application to groups and groupoids represents a significant advance in our understanding of these structures. Every diffeological group is, by its nature, a topological group when endowed with its D-topology, but it carries a much finer structure. It is this finer, smooth structure that allows us to distinguish and analyze objects that may be simple, or even trivial, from a purely topological viewpoint, yet are geometrically complex.

The quintessential example, and the one that motivated this entire research program, is the irrational torus. As a topological group, it is trivial, inheriting the coarse topology. As a diffeological group, however, it possesses a rich, non-trivial geometry, a non-trivial homotopy, and a fascinating group of diffeomorphisms whose structure depends on the arithmetic properties of the torus's slope. This demonstrates that the diffeological perspective reveals geometric and topological information that is invisible to topology alone.

The prime objects of ordinary differential geometry are manifolds, which are by definition locally homogeneous spaces. Their geometry, in the spirit of Klein's Erlangen Program, is captured by the action of global transformation groups. Diffeology, however, introduces a universe of spaces that are not locally homogeneous. These spaces possess intrinsic singularities, where the local geometric structure changes from one point to another. For such spaces, a single global group is no longer sufficient to capture the full richness of the geometry. A finer, more localized tool is required—one that can act differently at different points. This is precisely the role of the groupoid, the natural generalization of a group to an object with many units. While the use of groupoids in

*Date*: August 24, 2025.

*Key words and phrases.* Diffeology, Diffeological Groupoid, Stratification, Klein Groupoid, Orbifolds, Quasifolds, Noncommutative Geometry, Geometric Quantization.

The author thanks the Hebrew University of Jerusalem, Israel, for his continuous academic support. He is also grateful for the stimulating discussions and assistance provided by the AI assistant Gemini (Google).

diffeology is multifaceted—serving to define fiber bundles and other structures that refine classical constructions—an essential and indispensable role is to capture the local nature of singularities. This is the specific task of the **Klein groupoid**, the groupoid of local diffeomorphisms of a space.

The axiomatic framework of diffeology was first proposed by Jean-Marie Souriau in the early 1980s, initially tailored for the study of infinite-dimensional Lie groups [Sou80]. The focus of the theory shifted unexpectedly following a "crash test" of these axioms on the irrational torus [DIZ83]. It was in the process of building the necessary machinery to explain the non-trivial structure of this singular space that the central concept of this paper was elaborated [PIZ85]. The most economic and significative way to define a diffeological fiber bundle was through the concept of a *diffeological groupoid*.

This single idea proved to be the key. It not only solved the immediate problem of the irrational torus but also provided a powerful, unifying tool with broad applications. The diffeological groupoid became the natural language for constructing the universal covering of any diffeological space; for analyzing the intricate stratifications of orbifolds and orbit spaces; for building a crucial bridge to Alain Connes' noncommutative geometry; and, most recently, for forging a new approach to geometric quantization.

This paper is structured as a journey through these applications, illustrating the progressive realization of the groupoid's power.

- We will begin by revisiting the foundational work from 1985, explaining how the challenge of the irrational torus led to the definition of diffeological groupoids and, through them, a robust theory of fibrations and coverings.
- We will then explore its use in structural analysis, detailing how groupoids provide the right lens through which to view the stratification of orbifolds and the diffeology of orbit spaces.
- Subsequently, we will recount the pivotal role of the diffeological groupoid in forging a concrete link with noncommutative geometry, first for orbifolds and then for the more general class of quasifolds.
- A central contribution of this paper will be to make explicit a concept that has remained implicit in previous work: the *"Klein groupoid"* of a diffeological space, whose transitivity components define a fundamental, intrinsic stratification.
- Finally, to demonstrate the concept's continuing vitality, we will present its most recent application in the ongoing work on geometric quantization.

It is intended as a cohesive, expository narrative that not only recounts past results but also synthesizes them, highlighting the remarkable and consistent power of a single idea—the diffeological groupoid—to solve problems, unify concepts, and open new avenues of inquiry across geometry and mathematical physics.

## I. Foundations: The Diffeological Groupoid (1985)

The axiomatic of diffeologies, as first proposed by Jean-Marie Souriau, was not initially conceived to study singular spaces. Its purpose was to provide a rigorous differential framework for the most regular of objects: infinite-dimensional Lie groups, such as the group of symplectomorphisms or quantomorphisms [Sou80]. The focus of the theory shifted unexpectedly following a "crash test" of these new axioms on the quotient space of the irrational torus [DIZ83]. This space, which was topologically trivial, was expected to have a trivial diffeology. The surprising discovery that it possessed a rich and non-trivial smooth structure was the primary impetus for the subsequent research program. The immediate goal became to explain the previously computed homotopy groups of the irrational torus not as a standalone curiosity, but as a natural consequence of a general theory. The clear path was through the long exact homotopy sequence of a fibration. However, a significant obstacle stood in the way: the projection $T^2 \to T_\alpha$ is not locally trivial and does not satisfy the classical condition. At the time, no adequate theory of fiber bundles existed for this new category of spaces. It was in the process of building this necessary machinery, the subject of the 1985 thesis [PIZ85], that the central concept of this paper was elaborated. The most economic and significant way to define a diffeological fiber bundle was through the concept of a *diffeological groupoid* —a groupoid internal to the category of diffeological spaces.

At the heart of this story is the diffeological groupoid, but it appears in several distinct guises, each tailored to a specific geometric purpose. We will encounter *the structure groupoid*, an analytical tool derived from an atlas; *the Poincaré groupoid*, which captures the homotopy structure of paths; *the Klein groupoid*, which reveals the intrinsic singularity structure of any space; and finally, *the prequantum groupoid*, which encodes the geometric phase of a physical system. Understanding the roles and relationships of these different actors is key to appreciating the unifying power of the groupoid concept.

I.1. **The Diffeological Groupoid.** The foundational step was to equip the algebraic structure of a groupoid with a compatible smooth structure. A **diffeological groupoid** was defined as a groupoid **K** for which the set of objects, Obj(**K**), and the set of morphisms, Mor(**K**), are both diffeological spaces, satisfying four conditions of smoothness for the groupoid operations:

(1) The multiplication $(f, g) \mapsto f \cdot g$ is smooth.
(2) The inversion $f \mapsto f^{-1}$ is smooth.
(3) The source and target maps, src and trg, are smooth.
(4) The identity injection $x \mapsto \mathbf{1}_x$ is an *induction*.

An induction is a smooth map that is a diffeomorphism onto its image equipped with the subset diffeology. This definition ensures that all the structural maps of the groupoid respect the underlying diffeologies. A crucial consequence is that the space of objects Obj(**K**) can be identified diffeologically with the subspace of identity morphisms, meaning a diffeological groupoid is uniquely characterized by the diffeological space of its morphisms and the groupoid operation.

I.2. **A New Definition for Fiber Bundles.** With this new object, a more profound and powerful definition of fiber bundles became possible. For any smooth surjection $\pi : T \to B$, one can define its **structure groupoid**, whose objects are the points of the base B and whose morphisms are the diffeomorphisms between the fibers, $\text{Mor}(b, b') = \text{Diff}(T_b, T_{b'})$. This groupoid, endowed with a natural functional diffeology, captures the "smooth structure" of the projection.

It is important to acknowledge that the foundational idea of associating a groupoid with a fiber bundle was established by Charles Ehresmann in the late 1950s [Ehr59]. In the classical setting of differential geometry, he demonstrated a fundamental equivalence between principal bundles over a manifold and transitive Lie groupoids. This premonitory approach, however, was naturally limited by the category of manifolds.

The diffeological framework significantly enriches this perspective by extending its applicability in two crucial directions. First, the isotropy groups of a diffeological groupoid are not restricted to being Lie groups; they can be any diffeological group. This includes singular groups, such as the irrational tori that appear as fibers in the generalized prequantum bundles of non-integral closed 2-forms [PIZ95]. Second, and perhaps more profoundly, the structure groups can be infinite-dimensional. For a diffeological fibration $\pi : T \to B$ with fiber F, the isotropy groups of its structure groupoid are the diffeomorphism groups Diff(F), which are typically vast, infinite-dimensional diffeological groups. The concept of a fibrating groupoid is thus a far more general and flexible tool, capable of defining and analyzing a much wider universe of geometric structures.

This led to the central definition: a smooth projection $\pi : T \to B$ is a **diffeological fibration** if its structure groupoid is **fibrating**—that is, if its characteristic map $\chi : \mathrm{Mor}(\mathbf{K}) \to B \times B$ is a *subduction*. A subduction is a smooth surjection such that plots on the target space admit local smooth lifts into the source space at every point of their domain. This condition means that for any pair of points $(b, b')$ in the base, any plot of such pairs can be locally lifted to a smooth family of diffeomorphisms between the corresponding fibers.

This sophisticated definition was shown to be equivalent to a more operative one: a projection is a diffeological fibration if and only if it is **locally trivial along the plots** of the base. This equivalence was a key foundational result, establishing that the abstract groupoid condition precisely captured the necessary local lifting properties. Crucially, under this definition, the projection onto the irrational torus $T_\alpha$ is a diffeological fibration, thereby solving the initial problem and bringing this singular object squarely into a rigorous geometric framework.

I.3. **Principal Fibrations from Group Actions.** Within the general theory of diffeological fibrations, a particularly important subclass is that of *principal fibrations*, which are central to gauge theories and classification problems. The diffeological framework provides a remarkably direct and powerful criterion for determining when the quotient by a group action forms such a fibration.

Let G be a diffeological group acting smoothly on a diffeological space X. The key to understanding the structure of the quotient space X/G lies in the geometry of the action itself, captured by the **action map** (or graph of the action):

$$F : X \times G \to X \times X, \quad \text{defined by} \quad F(x, g) = (x, g(x)).$$

This map encodes how the group elements move the points of the space. The central result is the following theorem [PIZ85, PIZ13]:

**Theorem.** *The projection $\pi : X \to X/G$ is a principal diffeological fibration with structure group* G *if and only if the action map* F *is an **induction**.*

This condition essentially requires that the action is free (a consequence of F being injective) and that the orbits are "well-separated" in a diffeological sense. When this holds, the action is called a *principal action*.

This theorem is a powerful tool. It provides a direct criterion based on the group action itself to certify that the quotient map is a principal fibration, without needing to construct the full structure groupoid of the projection and verify that it is fibrating. It thus forms an essential, albeit indirect, link between the theory of group actions and the theory of groupoids in diffeology, showing how the properties of the former directly imply the fibrating nature of the latter.

I.4. **Immediate Successes of the Theory.** The definition of a diffeological fibration via fibrating groupoids was immediately validated by its powerful consequences, which directly addressed the problems that motivated its creation [PIZ85]. This foundational work established two crucial results.

First, it was proven that the projection of any diffeological group onto a quotient by any of its subgroups is a diffeological fibration. This general theorem is of paramount importance, as it instantly provides a rigorous framework for a vast class of quotient spaces, including the irrational torus. This result provided the formal justification for treating the projection $T^2 \to T_\alpha$ as a fibration, moving it from a singular case to a prime example of a general principle.

Second, and most crucially, this definition was precisely what was needed to establish the long exact sequence of homotopy in the context of diffeology. The existence of this sequence provided the ultimate conceptual tool. By applying it to the now-certified fibration $T^2 \to T_\alpha$, it provided a rigorous, conceptual confirmation of the initial by-hand computation of the fundamental group, $\pi_1(T_\alpha) \simeq \mathbf{Z}^2$. Furthermore, it allowed for the first formal computation of the higher homotopy groups—which had to be defined as part of this foundational work—showing that they vanish for $k > 1$. This necessity of building the theory from the ground up is reflected in the title of the thesis itself. The discovery of the diffeological groupoid thus successfully bridged the gap, providing the missing theoretical machinery.

I.5. **The Universal Covering via the Poincaré Groupoid.** The power of the diffeological groupoid was not limited to fiber bundle theory. It also provided a universal and elegant method for constructing covering spaces. The **Poincaré groupoid** (or fundamental groupoid) of a diffeological space X, denoted $\Pi(X)$, has the points of X as its objects and the fixed-ends homotopy classes of paths as its morphisms. The concatenation of paths endows $\Pi(X)$ with the structure of a diffeological groupoid which is, moreover, fibrating.

This insight allows for a direct construction of the universal covering. For a connected space X with a chosen basepoint $x_0$, the **universal covering space** $\tilde{X}$ is defined as the space of morphisms of the Poincaré groupoid starting at $x_0$:

$$\tilde{X} = \mathrm{Mor}_{\Pi(X)}(X, x_0, \star)$$

The projection $\pi : \tilde{X} \to X$ is given by the target map, $\pi([\gamma]) = \gamma(1)$. This construction yields a principal fiber bundle with the fundamental group $\pi_1(X, x_0)$ as its structure group. The main theorem of this construction asserts that for any connected diffeological space X, this procedure yields a unique (up to isomorphism) simply connected covering space.

This result is noteworthy for its remarkable generality, a feature that starkly contrasts with the classical topological setting. In topology, the existence of a universal covering requires the base space to satisfy local niceness conditions, such as being semilocally simply connected. The diffeological framework bypasses these requirements entirely. Because the Poincaré groupoid can be constructed for any diffeological space, the main theorem asserts that *every* connected diffeological space—regardless of whether it is singular, infinite-dimensional, or fails to satisfy classical local conditions—admits a unique, simply connected universal covering. This demonstrates the remarkable capacity of the diffeological groupoid framework to generalize fundamental topological constructions to a much broader class of spaces than classical manifolds.

II. Groupoids for Analyzing Structure: Stratification

The same conceptual machinery that proved so effective for defining global structures like fiber bundles also provides the essential tools for analyzing the local, layered structure of singular spaces. Many important diffeological spaces, such as manifolds with corners or orbifolds, are not homogeneous but are instead "modeled" on different local geometries from point to point [PIZ13, Chap. 4]. The intrinsic way to understand this layered structure is through the action of local diffeomorphisms, an action best formalized by a groupoid.

II.1. **The Structure Groupoid of an Orbifold.** To analyze the intrinsic structure of an orbifold X, we associate a groupoid to a choice of atlas. The construction proceeds in several steps [IZL18]:

(1) An **atlas** $\mathscr{A}$ for X is a collection of charts $f : \mathbf{R}^n/\Gamma \supset U \to X$ that cover X.
(2) From this, we form the **strict generating family** $\mathscr{F}$, consisting of the plots $F = f \circ \text{class} : \tilde{U} \to X$, where $\tilde{U} = \text{class}^{-1}(U) \subset \mathbf{R}^n$.
(3) The **nebula** $\mathscr{N}$ is the disjoint union of the domains of plots: $\mathscr{N} = \coprod_{F \in \mathscr{F}} \text{dom}(F)$. The nebula is a manifold, and the **evaluation map** $\text{ev} : \mathscr{N} \to X$, given by $\text{ev}(F, r) = F(r)$, is a subduction presenting X as a quotient of a manifold.

The **structure groupoid G** associated with the atlas $\mathscr{A}$ is then defined. Its objects are the points of the nebula, $\text{Obj}(\mathbf{G}) = \mathscr{N}$. Its morphisms are the germs of local diffeomorphisms of the nebula that are "invisible" at the level of the orbifold; that is, they project to the identity via the evaluation map:

$$\text{Mor}(\mathbf{G}) = \{\text{germ}(\Phi)_v \mid \Phi \in \text{Diff}_{\text{loc}}(\mathscr{N}) \text{ and } \text{ev} \circ \Phi = \text{ev} \restriction \text{dom}(\Phi)\}.$$

This groupoid captures the internal symmetries of the atlas presentation. The crucial link between this abstract construction and the geometry of the orbifold is the following fundamental lifting property [IKZ10]:

**Theorem (Lifting of Local Diffeomorphisms).** *Any local diffeomorphism of an orbifold lifts locally to a local diffeomorphism of the nebula.*

From this key result, everything else follows. It implies that the local symmetries of the orbifold are perfectly encoded by the arrows of the structure groupoid. The transitivity components of **G** are precisely the fibers of the evaluation map, meaning the space of components is diffeomorphic to the orbifold X itself. This leads directly to the classification of the intrinsic Klein strata of the orbifold: two points belong to the same Klein stratum if and only if their local isotropy groups within the structure groupoid are conjugate. This analysis culminates in the main theorem on the structure of orbifolds [GIZ23]:

**Theorem.** *The Klein stratification of a diffeological orbifold (with a locally finite atlas) is a standard stratification. Its strata are locally closed manifolds, and the partition is locally finite.*

It is worth noting how this geometric construction relates to other groupoid-based approaches to orbifolds, such as those of the Moerdijk school. In that framework, one often starts with a more abstract "orbifold groupoid," which may not be effective. It is a known result that the space of transitivity components of such a groupoid can be endowed with the structure of a diffeological orbifold. The structure groupoid **G** constructed here from the diffeology of the space is precisely the *effective* part of that starting object. This demonstrates that the diffeological approach is perfectly adequate for describing Satake's original V-manifolds, and that the equivalence class of the structure groupoid **G** emerges naturally as a fundamental invariant of the orbifold's geometry.

II.2. **The Diffeology of Orbit Spaces.** This analytical perspective extends to the more general case of orbit spaces M/G, where M is a manifold and G is a compact Lie group. The intrinsic geometry of the quotient space M/G is captured by its **Klein stratification**, the partition of the space into orbits of its pseudogroup of local diffeomorphisms. This stratification is fundamental because it decomposes the space into strata of equivalent singularity type, since points in the same stratum are, by definition, locally diffeomorphically equivalent.

The crucial insight of the diffeological approach is that this intrinsic partition is not merely a formal decomposition but is itself a *geometric stratification*: it is generated by the action of a groupoid. This object, which we will formally define in Section IV as the **Klein groupoid** of the space, has the local diffeomorphisms as its arrows, and its transitivity components are precisely the Klein strata.

The classical approach, by contrast, studies the projection of the *orbit-type stratification* from M to the quotient. This classical stratification on M groups points whose stabilizers are conjugate. A deeper analysis using diffeological tools reveals that this extrinsic picture is incomplete [GIZ25]. The intrinsic Klein stratification of the quotient does not always align with the projection of the orbit-type stratification. The correct correspondence is found by considering a finer partition on M: the **isostabilizer decomposition**. This partition groups points not just by the conjugacy class of their stabilizer, but by the stabilizer subgroup itself.

It is this decomposition that maps coherently to the intrinsic geometry of the quotient. The main result of this analysis is that the projection M → M/G induces a surjective map (not *a priori* injective as examples show) from the isostabilizer decomposition of M to the Klein stratification of M/G. The groupoid of local diffeomorphisms of the quotient is thus a more refined instrument than the classical orbit-type structure on the source manifold. It probes the *intrinsic* geometry, revealing the true singularity structure and establishing the correct correspondence between the group action on M and the diffeological geometry of the resulting orbit space M/G.

Thus, while the structure groupoid of an atlas is the key to understanding the presentation of an orbifold, it is the Klein groupoid that reveals the intrinsic geometry of more general orbit spaces.

### III. A Bridge to Noncommutative Geometry

One of the most fruitful applications of the diffeological groupoid has been in building a concrete, structural bridge to Alain Connes' Noncommutative Geometry. The core philosophy of this program is to study "singular" geometric spaces by associating them with non-commutative algebras, typically C*-algebras. The diffeological groupoid, constructed directly from the geometry of the space, provides the perfect intermediate object to make this association rigorous and functorial.

III.1. **The Orbifold Blueprint.** The strategy for connecting diffeology to noncommutative geometry was first elaborated and proven in the case of orbifolds [IZL18]. As established in the previous section, any atlas of an orbifold X gives rise to a **structure groupoid G**. This groupoid is not just an abstract entity; it is a well-behaved diffeological groupoid which is both *étale* and *Hausdorff*.

These properties are precisely the technical requirements needed to apply the standard construction, developed by Jean Renault in [Ren80], for associating a C*-algebra to a groupoid. The construction proceeds by considering the space of compactly supported, continuous complex-valued functions on the space of morphisms, Mor(**G**). This space is endowed with an algebraic structure where the groupoid multiplication defines a convolution product, and the groupoid inversion defines an involution. The completion of this convolution algebra with respect to a suitable norm yields the **C\*-algebra of the groupoid**, denoted C*(**G**).

This provided the essential blueprint: a direct path from a geometric object (an orbifold, via its atlas) to an algebraic one (a C*-algebra). It realized the goal of replacing a singular space with a well-behaved, albeit non-commutative, algebra.

III.2. **Generalization to Quasifolds and Morita Equivalence.** The full power of this connection is most evident in the generalization to **quasifolds** [IZP21]. Quasifolds are a broader class of diffeological spaces, locally modeled on quotients $\mathbf{R}^n/\Gamma$ where $\Gamma$ is a countable (but possibly infinite) subgroup of the affine group. This class includes the most interesting examples, such as the irrational torus $T_\alpha$, and was originally explored by Elisa Prato inspired by the study of quasicrystals [Pra99].

The primary challenge in extending the construction from orbifolds to quasifolds was that the local structure groups $\Gamma$ are no longer finite. Their orbits can be dense, and the proofs required substantial revision. The key technical insight that made this generalization possible was the use of Baire's theorem, which allowed arguments to be extended from the finite to the countable case.

This generalization led to a crucial question. The construction of the C*-algebra depends on the choice of an atlas for the quasifold. A different atlas produces a different nebula, a different structure groupoid **G**′, and thus a different C*-algebra C*(**G**′). For the construction to be meaningful, these different algebras must represent the same underlying "noncommutative space." This is captured by the central result of the theory:

**Theorem.** *Let* **G** *and* **G**′ *be the structure groupoids associated with two different atlases of the same quasifold* X. *Then the corresponding C\*-algebras,* C*(**G**) *and* C*(**G**′)*, are **strongly Morita equivalent***.

Morita equivalence is the accepted notion of isomorphism for C*-algebras in noncommutative geometry. This theorem thus establishes that the construction yields a true invariant of the underlying space. It associates a unique Morita equivalence class of C*-algebras to the diffeomorphic class of the quasifold. This result solidifies the bridge between the two theories, confirming that the geometric information encoded in the diffeology of a quasifold is faithfully preserved in the algebraic structure of its associated noncommutative algebra.

This result solidifies the bridge between the two theories, confirming that the geometric information encoded in the diffeology of a quasifold is faithfully preserved in the algebraic structure of its associated noncommutative algebra. This naturally opens avenues for further research aimed at deepening the dictionary between geometry and algebra. For instance, a compelling question, prompted by discussions with Jiawen Zhang of Fudan University, is to find the geometric interpretation within the quasifold X of the ideal structure of its C*-algebra. For groupoid C*-algebras, this structure is known to correspond to the open, invariant subsets of the groupoid's space of objects [Ren80]. Understanding how these subsets, which live on the nebula, project to a meaningful geometric structure on the quotient space X is a subtle and intriguing problem for future investigation.

## IV. The Klein Groupoid: An Explicit Framework for Intrinsic Stratification

In the classical geometry of manifolds, one can adopt two complementary perspectives, a distinction often discussed with Souriau. The first is a **passive view**, where the manifold is understood through an atlas of charts—a collection of external parameterizations. The second is an **active view**, where the manifold is understood as a "theater" for the action of its group of diffeomorphisms. In diffeology, this distinction becomes even more acute and powerful. The passive view is embodied by the set of all plots, which forms the complete, maximal atlas for the space. The active view, however, requires a more refined tool than the global diffeomorphism group, especially for spaces that are not locally homogeneous.

For these singular spaces, where the geometry changes from point to point, a global group of transformations is too coarse an instrument. A finer tool is needed, one that captures the local nature of symmetry. This is the essential role of the groupoid. The **Klein groupoid** of a diffeological space is the ultimate realization of this active viewpoint. It is the groupoid of all local symmetries, and it is through its structure that the true intrinsic geometry of the space reveals itself.

The previous sections have demonstrated the power of groupoids tailored to specific structures: the structure groupoid of a fibration, the structure groupoid of an orbifold atlas, or the groupoid of local diffeomorphisms of a quotient space. We now introduce a canonical and fundamental groupoid that can be associated with *any* diffeological space, providing a universal tool for understanding its intrinsic geometric structure. This concept, which has been implicit in much of the preceding analysis, we now make explicit.

IV.1. **Formal Definition.** The local geometry of a diffeological space is characterized by its pseudogroup of local diffeomorphisms. It is natural to formalize this structure into a groupoid.

**Definition (The Klein Groupoid).** *Let* X *be any diffeological space. The **Klein groupoid** of* X, *denoted* $K\ell(X)$, *is the diffeological groupoid whose objects are the points of* X *and whose morphisms are the germs of local diffeomorphisms of* X.

- $\text{Obj}(K\ell(X)) = X$.
- $\text{Mor}(K\ell(X)) = \{\text{germ}(\varphi)_x \mid \varphi \in \text{Diff}_{\text{loc}}(X), x \in \text{dom}(\varphi)\}$.

*The source and target maps are given by* $\text{src}(\text{germ}(\varphi)_x) = x$ *and* $\text{trg}(\text{germ}(\varphi)_x) = \varphi(x)$.

The smooth structure on this groupoid is not trivial. The space of morphisms is endowed with a specific **functional diffeology** that was formally constructed for this purpose in [IZL18]. The construction is a two-step process: first, a functional diffeology is defined on the space of local smooth maps itself, a structure that carefully handles the variable domains of definition. The space of morphisms $\text{Mor}(K\ell(X))$ is then endowed with the quotient diffeology induced by the "germification" map, which takes a local diffeomorphism and a point in its domain to the corresponding germ.

This groupoid is the ultimate embodiment of Felix Klein's Erlangen Program in the context of diffeology. It captures the full local symmetry of the space, where "geometry" is defined by the properties invariant under local diffeomorphisms.

IV.2. **The Klein Stratification and its Invariants.** The action of the Klein groupoid on its space of objects naturally partitions the space into its orbits.

**Definition (The Klein Stratification).** *The Klein stratification of a diffeological space* X *is the partition of* X *into the transitivity components (orbits) of its Klein groupoid* $K\ell(X)$.

This stratification immediately captures the most fundamental numerical invariant of a diffeological space. Since local diffeomorphisms preserve the diffeological dimension, we have the following foundational result:

**Theorem.** *The diffeological dimension map is constant on each Klein stratum.*

The dimension map is thus the first and most basic invariant revealed by the action of the Klein groupoid. The partition into orbits is also a true stratification in the topological sense. Because local diffeomorphisms are also homeomorphisms for the D-topology, this partition automatically satisfies the frontier condition, making it a *basic stratification*. The Klein strata are precisely the sets of points that are locally diffeomorphically equivalent, thus partitioning the space into its distinct singularity types, each with a well-defined dimension.

IV.3. **Unifying Power.** The true power of this concept lies in its ability to provide a unified perspective on the structures discussed in the previous sections.

For a smooth manifold, the Klein groupoid is transitive. The action of local diffeomorphisms can connect any point to any other point, and thus the Klein stratification consists of a single stratum: the manifold itself. This reflects the principle of homogeneity for manifolds.

For more complex spaces, particularly those constructed as quotients, the Klein groupoid provides the definitive language for their intrinsic geometry. It replaces formal, extrinsic structures with a geometric one generated by a groupoid action. Consider the case of an orbit space M/G. As we saw, the intrinsic Klein stratification of the quotient reveals the true singularity types. This intrinsic structure "encrypts" the essential geometric information of the group action on M, with the isostabilizer decomposition being the precise structure on the source that maps coherently to it.

This principle is not reserved for a small class of examples; it is a fundamental tool in diffeology for understanding the geometry of a space. In an analogous manner, if we consider a space whose diffeology is the pushforward of a projection map—for example, the "shadow" of a surface Σ projected onto a plane—the Klein stratification of the resulting space will necessarily encryp until some level the geometric information about the singularities of the projection map.[1]

In all these cases, the Klein groupoid acts as the ultimate probe of the local structure. It does not depend on how the space was constructed, whether by a group action, a projection, or some other means. It operates directly on the final space and reveals its inherent geometric hierarchy. In this sense, the Klein groupoid and its stratification are not just tools for studying the geometry; they *are* the intrinsic local geometry of the space.

---

[1]This is certainly a path to follow and clarify.

IV.4. **Subgroupoids and Geometric Structures.** The Klein groupoid $K\ell(X)$ captures the full local symmetry of a diffeological space. A powerful feature of this framework is that by restricting to subgroupoids, we can analyze additional geometric structures on X. If X is equipped with a geometric structure 'S' (such as a Riemannian metric [PIZ25a, p.224] or a parasymplectic form),[2] we can consider the subgroupoid of the Klein groupoid whose morphisms are the germs of local diffeomorphisms that preserve 'S'. The transitivity components of this subgroupoid then classify the points of X according to their local geometry with respect to 'S'.

This approach is particularly fruitful in symplectic diffeology. Let $(X, \omega)$ be a parasymplectic space. We can define the **symplectic Klein groupoid**, $K\ell(X, \omega)$, as the subgroupoid of $K\ell(X)$ consisting of germs of local automorphisms of $(X, \omega)$.

This leads to a purely geometric and diffeological definition of a presymplectic space. In the theory of symplectic diffeology [PIZ10], a parasymplectic space $(X, \omega)$ is defined as **presymplectic** if its symplectic Klein groupoid $K\ell(X, \omega)$ is transitive. For manifolds, this condition of local homogeneity is equivalent to the classical definition of a presymplectic form having constant rank. However, this new definition is a more fundamental generalization. It replaces the classical criterion—a condition on the kernel of the form, which relies on the tangent bundle—with an intrinsic, geometric one. This avoids the often problematic and unnecessary complexity of defining tangent spaces for general diffeological spaces, providing a robust definition that is native to the category.

This perspective connects naturally to the theory of moment maps. The universal moment map, $\mu_\omega$, is a map from the space X to the space of momenta of the *global* group of automorphisms, $\text{Diff}(X, \omega)$. This global object motivates a compelling avenue for future research: the question of a "local moment map" living naturally on the symplectic Klein groupoid itself. Its construction is a central problem in symplectic diffeology, as it would provide the key to extending the fundamental theorem that identifies the characteristics of a homogeneous presymplectic space with the level sets of the universal moment map [PIZ13, §9.26]. Such a generalization is crucial for infinite-dimensional models, where the classical kernel-based definition of characteristics is problematic.

IV.5. **A Functorial Viewpoint: The Complete Picture.** The relationship between the geometry of the nebula $\mathcal{N}$ (a manifold) and the quasifold X (a singular space) can be framed in a powerful, functorial way that provides the complete picture. We have two fundamental geometric objects: the Klein groupoid of the nebula, $K\ell(\mathcal{N})$, and the Klein groupoid of the quasifold, $K\ell(X)$.

Consider the subgroupoid of $K\ell(\mathcal{N})$, let us call it $\mathbf{G}_{\text{ev}}$, whose morphisms are the germs of local diffeomorphisms $\Phi$ of $\mathcal{N}$ that project to local diffeomorphisms of X. That is, a morphism $\text{germ}(\Phi)_\nu$ belongs to $\mathbf{G}_{\text{ev}}$ if there exists a local diffeomorphism $\varphi$ of X such that the following diagram commutes locally:

$$\begin{array}{ccc} \mathcal{N} & \xrightarrow{\Phi} & \mathcal{N} \\ \text{ev} \downarrow & & \downarrow \text{ev} \\ X & \xrightarrow{\varphi} & X \end{array}$$

---

[2]We call parasymplectic any closed 2-form.

This construction naturally defines a functor of groupoids, $\text{ev}_* : \mathbf{G}_{\text{ev}} \to \text{K}\ell(X)$. The fundamental Lifting Theorem for orbifolds and quasifolds [IKZ10, IZP21] states precisely that this functor is surjective on morphisms. Since the evaluation map ev is surjective by construction, the functor is also surjective on objects.

The Lifting Theorem implies that the functor $\text{ev}_*$ is **surjective on both objects and morphisms**. This is a very strong condition, much stronger than a mere equivalence of categories, and it establishes that the Klein groupoid of the quasifold is, in a precise sense, a quotient of the lifted groupoid from the nebula.

This provides a complete, hierarchical picture. The intrinsic geometry of the manifold $\mathcal{N}$ projects onto the intrinsic geometry of the singular space X. The **structure groupoid G** of the atlas, which we used to build the C*-algebra, finds its ultimate conceptual place here: it is precisely the kernel of this functor $\text{ev}_*$. It is the "gauge" subgroupoid consisting of all the internal symmetries of the atlas presentation that project to the identity (the units) in the Klein groupoid of the quasifold.

This functorial picture can be generalized to any subduction, raising a deep structural question: what property of a subduction ensures that the associated functor is surjective on morphisms? This "Lifting Property for Local Diffeomorphisms" is not true for all subductions. It holds for diffeological fibrations, where the uniformity of the fibers allows for lifting. More subtly, it holds for the subduction $\text{ev} : \mathcal{N} \to X$ for quasifolds, which is not a global fibration but rather a diffeological foliation whose transverse structure is locally constant, modeled on quotients $\mathbf{R}^n/\Gamma$. This suggests that the class of subductions that admit such a lifting are those that possess a certain "local fibration" or "uniform foliation" structure. Identifying and classifying this family of subductions is a compelling direction for future research.

In this sense, the Klein groupoid and its stratification are not just tools for studying the geometry; they *are* the intrinsic local geometry of the space. It is this power to capture the most fundamental geometric essence of a space that allows the groupoid concept to be applied to the most profound problems, including those at the interface of geometry and physics, as we shall now see.

## V. A New Frontier: Geometric Quantization by Paths

The journey of the diffeological groupoid culminates in its most recent application: forging a new, unified path in geometric quantization. This approach is deeply connected to the spirit of Feynman's path integral formulation of quantum mechanics, which posits that the quantum amplitude is a sum over all possible histories (paths) of a system. Our construction provides a geometric framework for this idea, building the prequantum structure directly from the space of paths of the classical system.

V.1. **The Prequantum Groupoid: A Quantization à la Poincaré.** Classical geometric quantization faces a well-known topological obstruction: a "prequantum line bundle" only exists if the symplectic form is integral. In a forthcoming work [PIZ25b], a new prequantum object is constructed for any connected and simply connected parasymplectic space $(X, \omega)$ (a diffeological space equipped with a closed 2-form) that bypasses this obstruction entirely. The new structure is not a bundle, but a **prequantum groupoid**, $\mathbf{T}_\omega$, built as a diffeological quotient of the space of paths on X.

This construction is best understood as a "quantization" à la Poincaré.

- The **Poincaré groupoid** Π(X) is a quotient of the space of paths where two paths are equivalent if they are homotopic. Its isotropy group at a point $x$ is the fundamental group $\pi_1(X, x)$, which captures the *homotopic phase*.
- The **prequantum groupoid** $\mathbf{T}_\omega$ is a quotient of the space of paths where two paths are equivalent if the integral of $\omega$ over any surface connecting them is a "period" of $\omega$. Its isotropy group at a point $x$ is the torus of periods $\mathrm{T}_\omega$, which captures the *geometric phase*.

The prequantum groupoid thus replaces the homotopic equivalence with a geometric equivalence dictated by the parasymplectic form.

V.2. **The Geometry of the Prequantum Groupoid.** The space of morphisms of the prequantum groupoid, $\mathscr{Y}_\omega = \mathrm{Mor}(\mathbf{T}_\omega)$, has a rich geometric structure. The "ends" map, $\mathrm{ends} : \mathscr{Y}_\omega \to \mathrm{X} \times \mathrm{X}$, is a diffeological fibration. The fiber over a pair of points $(x, x')$ is the space of all equivalence classes of paths from $x$ to $x'$. This fiber is not a group, but it is a *principal homogeneous space* (a torsor) for the torus of periods $\mathrm{T}_\omega$. This means that while there is no canonical "identity" path, the group $\mathrm{T}_\omega$ acts freely and transitively on the set of paths connecting any two points.

It is important to note that this fibration is not, in general, a principal bundle. A principal bundle requires a global group action on the total space, whereas here the action is "fiber-wise" by the isotropy groups. This fibration of $\mathrm{T}_\omega$-torsors over $\mathrm{X} \times \mathrm{X}$ is the fundamental geometric object of the prequantization.

V.3. **Advantages and Future Directions.** This path-based construction has two profound advantages over the classical approach. First, the prequantum groupoid **always exists** for any parasymplectic space with discrete periods.[3] Second, it **preserves all symmetries**: the full group of automorphisms of the parasymplectic structure, $\mathrm{Diff}(X, \omega)$, acts isomorphically as a group of automorphisms on the prequantum groupoid structure.

This framework opens several avenues for a "true" quantization procedure. The natural candidate for the space of quantum states, or wave functions, is the space of *multiplicative functions* on the groupoid—smooth functions $\psi : \mathscr{Y} \to \mathbf{C}$ such that $\psi(y_1 \cdot y_2) = \psi(y_1)\psi(y_2)$.[4] The group of automorphisms $\mathrm{Diff}(X, \omega)$ acts naturally on this space, providing the quantum dynamics.

The power of this diffeological construction is its vast applicability. It applies uniformly to any diffeological space, be it a standard manifold, an orbifold, a quasifold, or an infinite-dimensional space. A current investigation, for example, involves applying this framework to the quantization of the deuteron, modeled as a system on the orbifold $\Sigma^4 = (S^2 \times S^2)/\mathfrak{S}_2$. This work demonstrates the continuing vitality of the diffeological groupoid, showing that this concept, born from a need to understand singular fibrations, has evolved into a fundamental tool capable of providing new solutions and deeper insights into the most current problems in geometry and mathematical physics.

By building the quantum phase structure directly from an equivalence relation on paths, this approach offers a compelling geometric realization of the physical intuition behind the path integral.

---

[3] Here, discrete periods means a diffeologically discrete subgroup of $\mathbf{R}$, that is, any strict subgroup $\mathrm{P}_\omega \subsetneq \mathbf{R}$. This is the case in particular of all Hausdorff second countable manifolds.

[4] In the ususal case of prequantization, when the torus of periods is $\mathrm{U}(1)$.

**Note on the General Case.** The construction presented here for simply connected spaces serves as the foundation for the general case. In a subsequent paper in progress, the theory is extended to arbitrary connected diffeological spaces by lifting the structure to the universal covering space and then taking a quotient by the action of the fundamental group. This general construction reveals a richer isotropy structure, where the periods of ω are some extension, related to the fundamental group $\pi_1(X)$, of the periods of the lifted 2-form.

## Conclusion

This expository paper has recounted a thirty-year journey, tracing the evolution of a single, powerful concept: the diffeological groupoid. Born from a specific challenge—the need to make sense of the irrational torus—it has grown into a fundamental tool with applications across a wide swath of modern geometry.

We have seen the groupoid play several distinct but interconnected roles. As a **definitional tool**, the concept of a *fibrating groupoid* provided the key to a robust theory of fiber bundles and coverings in diffeology, one capable of rigorously explaining the homotopy of singular spaces and establishing the long exact sequence.

The concept then evolved into a powerful **analytical instrument**. The *structure groupoid* of an atlas allowed for a deep analysis of the internal symmetries of orbifolds and quasifolds. More fundamentally, the *Klein groupoid* of local diffeomorphisms was shown to provide a universal language for the intrinsic geometry of any diffeological space, with its transitivity components—the Klein strata—partitioning the space into loci of constant diffeological dimension and equivalent singularity type.

This analytical power, in turn, allowed the groupoid to serve as a **bridge between theories**. It forged a concrete, functorial link between the geometric world of singular diffeologies and the algebraic world of noncommutative geometry, associating a unique Morita equivalence class of C*-algebras to each quasifold. Most recently, it has offered a new bridge between classical and quantum mechanics. The *prequantum groupoid*, constructed as a "quantization" in the mode of the Poincaré groupoid, provides a geometric realization of the path-integral philosophy, resolving long-standing issues in geometric quantization.

From a definitional tool to an analytical one, from a bridge between fields to a unifying principle, the diffeological groupoid has proven to be an exceptionally durable and fertile concept. Its ability to solve problems, clarify structures, and open new avenues of inquiry across a wide range of topics in geometry and mathematical physics is a testament to its fundamental nature. The journey is far from over, and the full potential of this powerful idea —–to unify, to clarify, and to solve—– is surely yet to be exhausted

PATRICK IGLESIAS-ZEMMOUR, EINSTEIN INSTITUTE OF MATHEMATICS, THE HEBREW UNIVERSITY OF JERUSALEM, EDMOND J. SAFRA CAMPUS, GIVAT RAM, 9190401 JERUSALEM, ISRAEL.

*Email address*: piz@math.huji.ac.il

*URL*: http://math.huji.ac.il/~piz